# ELECTRIFYING THE URBAN TAXI FLEET: A DATA-DRIVEN APPROACH

Yikang Li[a], Huimiao Chen[b], Yanni Yang[c], Ziyang Guo[d], Fang He[a,*]

[a]Department of Industrial Engineering, Tsinghua University, Beijing, 100084, PR China

[b]Department of Electrical Engineering, Tsinghua University, Beijing 100084, PR China

[c]Information School, Capital University of Economics and Business, Beijing 100070, PR China

[d]Sparkzone Institute, Beijing 100084, PR China

December 10, 2017

Abstract: This paper is devoted to proposing a data-driven approach for electrifying the urban taxi fleet. Specifically, based on the gathered real-time vehicle trajectory data of 39,053 taxis in Beijing, we conduct time-series simulations to derive insights on both the configuration of electric taxi fleet and dispatching strategies. The proposed simulation framework accurately models the electric vehicle charging behavior from the aspects of time window, charging demand and availability of unoccupied charges, and further incorporates a centralized and intelligent fleet dispatching platform, which is capable of handling taxi service requests and arranging electric taxis' recharging in real time. To address the impacts of the limited driving range and long battery-recharging time on the electrified fleet's operations efficiency, the dispatching platform integrates the information of customers, taxi drivers and charging stations, and adopts a rule-based approach to achieve multiple objectives, including reducing taxi customers' average waiting time, increasing the taxi demand fill rate and guaranteeing the equity of the income among taxi drivers. Although this study only examines one type of fleet in a specific city, the methodological framework is readily applicable to other cities and types of fleet with similar dataset available. Simulation results indicate that: i) after EV taxi number is greater than 12,000, station capacity is larger than 16, and the battery range is more than 200km, continuously increasing them brings marginally diminishing effect for electrifying the taxi fleet in Beijing; ii) choosing a sufficient battery range could not only enhance the efficiency of the system but also promote the fairness of drivers' income; iii) the dispatching strategy should dynamically fit the electric taxi fleet configuration as well as customer demand level; iv) the curve of charging demand shows a much more fluctuated pattern than the curve of customer demand.

---

* Corresponding author. Email: fanghe@tsinghua.edu.cn.

**1. INTRODUCTION**

Electric vehicles (EVs) receive rapidly increasing attention because of the advancement of battery technologies and the growing concern over environment issues. (e.g., Larminie and Lowry, 2003; He et al., 2013abc, 2014, 2015 and Gardner et al., 2013). Compared to conventional gasoline-powered vehicles, EVs hold several advantages, including increasing energy security, achieving better fuel economy and reducing tailpipe emissions, among others (US Department of Energy, 2013; Ghamami et al., 2015; Kontou et al., 2015). Many government agencies are promoting the deployment of EVs through providing incentive policies, such as offering purchase subsidies and deploying more public charging stations in urban areas.

Taxi fleet serves as an important component of mobility in the transportation system of any modern city. In recent years, multiple cities and metropolitan areas around the world have started electrifying the taxi fleet, i.e., adopting the taxi fleet composing of EVs (e.g., Sathaye 2013; Baptista et al., 2011; Gao and Kitirattragarn, 2008; Reynolds et al., 2011). Although the electric taxi fleet also enjoys the aforementioned advantages of EVs, its efficiency will be possibly affected by the limited driving range of onboard batteries and long battery-recharging time (Lu et al., 2012; Gacias and Meunier, 2015; Li et al., 2017; Yang et al., 2016). To tackle this challenge, a promising solution would be optimally designing EV battery range, deploying charging infrastructure and developing a real-time taxi fleet dispatching platform with the consideration of EVs' limited range and recharging plan.

To optimally deploy public charging infrastructure, various approaches have been proposed in the literature. The flow-capturing models optimize the charging station locations to maximize the amount of travelers passing by charging stations (e.g., Hodgson, 1990; Berman et al., 1992, 1995; Hodgson and Berman, 1997). The network-equilibrium approach captures the EV drivers' reactions to the deployment of public charging stations and other travelers' decisions, and optimizes the station location to maximize the social welfare (e.g., He et al., 2013a; 2015; Jiang et al., 2012; Jiang and Xie, 2014; Chen et al., 2016; Wang et al., 2016). However, both above approaches are based on some assumptions of EV drivers' behaviors, which remain to be verified by the real-world data. Recently, real-world driving profiles have been utilized to estimate EVs' public charging needs and then optimize the station locations (e.g., e.g., Dong et al., 2013; Andrews et al., 2012; Dong and Lin, 2012). Among them, some studies (e.g., Cai and Xu, 2013; Cai et al., 2014; Yang et al., 2017; Shahraki et al., 2015; Yang et al., 2016) optimally locate charging stations through applying the "big data" mining techniques on the real-time and large-scale trajectory data, and revealing the inherent heterogeneity of individual travel patterns. Note that the above studies assume that the travel profiles of taxi drivers remain unchanged after shifting from conventional gasoline-powered vehicles to EVs, which is a strong assumption because in reality, EV taxi drivers may totally change their travel pattern based on their specific battery states of charge and different charging station utilization levels. Furthermore, because of the limited battery range, it is highly possible that EV taxi drivers prefer to rely on the e-hailing or ride-sourcing smartphone-applications (such as Uber or Didi) to search for customers, rather than waste battery ranges on cruising on streets.

With regard to optimally dispatching the gasoline-powered taxi fleet, many studies have been published in the literature. One group of studies focuses on dividing the area into zones and applying rule-based approach to dispatch taxis to satisfy the service requests (e.g., Shrivastava et al., 1997; Alshamsi et al., 2009; Seow et al., 2010). Another group of studies considers that the taxi customers make reservation in advance, and models the taxi dispatching problem as a variant of the pick- up and delivery vehicle routing problem with time windows (e.g., Wang et al., 2009; Horn, 2002; Meng, et al, 2010). Compared to gasoline-powered taxi fleet, the number of studies on how to dispatch the electric taxi fleet is limited. Among the existing studies, Wang and Cheu (2013) only focus on the service requests with advance reservations, whereas Gacias and Meunier (2015) take into account the mixed service demands of both advance reservations and street hailing. Sathaye (2014) provides an analytic optimization framework for the design of electric taxi system. In order to reduce the taxi drivers' recharging time, Lu et al. (2012) propose a dispatching policy which comprehensively considers the taxi demand, the remaining charge of EVs' battery, and availability of charging station. However, none of the above studies reveals the impacts of the taxi fleet configuration, charging infrastructure deployment plan and demand level on the design of taxi dispatching strategies. Moreover, they also ignore examining the taxi drivers' income level, which is actually important for maintaining the urban taxi system's sustainability.

Inspired by the above studies, this paper is devoted to proposing a data-driven approach for electrifying the urban taxi fleet. In order to characterize the spatial and temporal distribution of customer demands, this paper gathers the real-time vehicle trajectory data of 39,053 taxis in Beijing in May of 2016, and extracts the customer origin-destination (O-D) pairs whose number reaches approximately eight million in one month. On the basis of this dataset, we conduct time-series simulations to derive insights on both the configuration of electric taxi fleet and dispatching strategies. Specifically, in the simulation, the electric vehicle recharging behaviors are accurately modeled from the aspects of time window, charging demand and availability of unoccupied charges. Moreover, a centralized fleet-dispatching platform, similar to Uber or Didi, is developed and considered to serve taxi service requests and arranging electric taxis' recharging in real time. In such context, taxi customers request rides via the mobile application, and upon receiving the ride requests, the centralized platform will match the electric taxis with the taxi customers. To address the impacts of the limited driving range and long battery-recharging time on the electrified fleet's operations efficiency, the dispatching platform integrates the information of customers, taxi drivers and charging stations, and adopts a rule-based approach to achieve multiple objectives, including reducing taxi customer average waiting time, increasing the taxi demand fill rate and guaranteeing the equity of the income among taxi drivers. Simulation results quantify the impacts of EV number, battery range, station capacity and dispatching-strategy parameters on demand fill rate, customer average waiting time and Gini coefficient of taxi driver income.

Compared to previous studies (e.g., Cai et al., 2014; Cai and Xu, 2013; Lu et al., 2012), the contribution of this paper lies in the following four aspects. Firstly, we consider that the number of chargers at each charging stations is limited. Therefore, EVs can charge batteries only if there are still unoccupied chargers left at stations, and otherwise they need to wait in the queue for next available chargers. Doing this makes our simulations capable of accurately modeling the real-time operations of public charging stations and reflecting the interactions of different EVs' charging behaviors. Secondly, we relax the unrealistic assumption of the unchanged travel pattern of taxi drivers after shifting to EVs, through developing a centralized fleet-dispatching platform. To facilitate the real-time operations, the platform adopts a rule-based approach and integrates the information of customers, taxi drivers and charging stations. Thirdly, the performance of electrified taxi fleet is examined from the aspects of taxi customer average waiting time, taxi demand fill rate and the equity of the income among taxi drivers. Note that the taxi-driver income-level distribution is not widely considered by the previous studies. However, we believe it is an important factor for the electric taxi fleet to consider, because electric taxi drivers' income could be potentially affected by the limited battery range and possibly long recharging time. Fourthly, we quantify the impacts of EV number, battery range, station capacity and dispatching-strategy parameters on the performances of electrified taxi fleet, simultaneously.

For the remainder of this paper, section 2 introduces the dataset and the method of deploying charging stations. Section 3 derives the time-series simulations and develops the real-time taxi fleet dispatching platform. In section 4, different simulation results are analyzed to derive the insights on the impacts of EV number, battery range, station capacity and dispatching-strategy parameters. Section 5 concludes the paper.

## 2. DATA AND CHARGING STATION DEPLOYMENT

Using Beijing as a case study, we extract the dynamic O-D taxi customer demands from the vehicle trajectory dataset of 39,053 taxis. Spatially, the O-D demands distribute over the entire city, and their number could reach approximately eight million in four weeks. We are interested in exploring how to electrify the urban taxi fleet to serve such set of O-D demands. Because of the limited driving range and long battery recharging time, the deployment of charging infrastructure plays a key role in the operations of the electrified taxi fleet. In this section, on the basis of this dataset, we will investigate on how to deploy public charging infrastructure.

## 2.1 Data Description and Preprocessing

To better characterize the spatial and temporal distribution of O-D taxi customer demands at a city level, we gather the real-time vehicle trajectory data of 39,053 taxis in Beijing from May 1st to May 27th in 2016, collected by smartphone and on-board device[1]. The dataset records the state (including location, speed and whether a taxi is in service) of individual taxi every 30 seconds. We first group the data by vehicle ID and sort the data in each group by time. To extract the real-time taxi customer demands from the vehicle trajectory set, we devote substantial efforts to cleaning the data. For instance, if a trip' origin or destination lacks the information of GPS coordinates, we approximate it using the nearest data point that is within three minutes from the trip's start or end time, respectively. We mark the trips, whose travelling time is less than two minutes, as invalid trips.

Figure 1 shows the histogram of the total number of taxi trips in Beijing within one day from our dataset. It can be observed that the total trip number varies from 220,000 to 340,000, with the mean and standard deviation values as 296,230 and 37,469, respectively. We randomly select 25% of the trips and plot each trip's origin on the map, as shown in Figure 2.

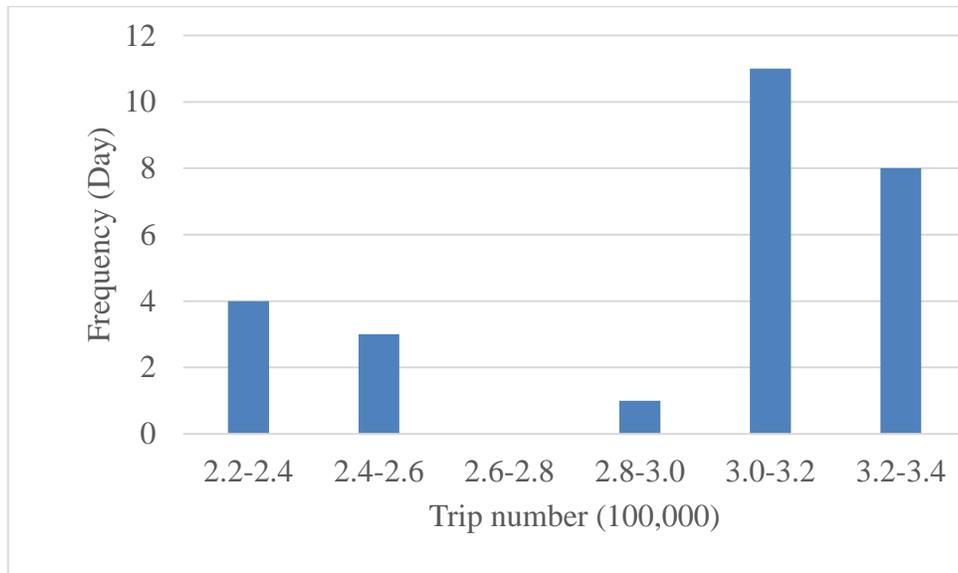

Figure1. Histogram of the total number of one-day taxi trips in the dataset

---

[1] The total number of taxis in Beijing is approximately 66,000 (Huo et al., 2012; Zheng et al., 2011).

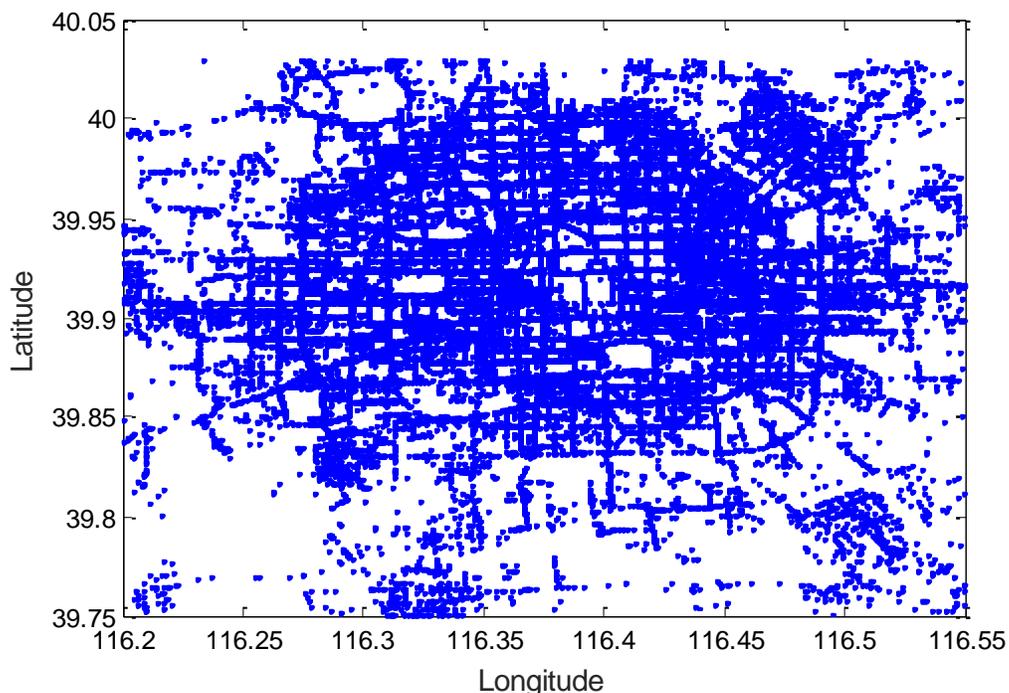

Figure2. The distribution of trips' origins

**2.2 Charging Station Deployment**

The distributions of trips' origins and destinations are significant factors for deploying charging stations, because electric taxis can only recharge at the intermediate stops between two adjacent trips. We apply K-means to cluster the locations of trips' origins, and then locate charging station to each cluster at its centroid[2]. Note that this locating method is in accordance with the suggestion by Cai et al. (2014) that the number of parking events serves as a good criterion to locate stations. For instance, if we plan to locate 100 public charging stations, we will apply the K-means method to partition the locations of all the trips' origins into 100 clusters and then locate a station at the centroid of each cluster. Figure 3 shows an example, where each color represents a particular cluster and each cross represents a charging station. Obviously, the stations are densely distributed in the urban area. The design of station capacity will be discussed in section 4.

---

[1] We observe that the distribution of taxi trips' origins and destinations are similar. Hence, we only utilize the distribution of origins to deploy charging stations. Furthermore, we do not require the station locations to sit in the existing gasoline stations in consideration of the fact that the existing gasoline stations do not necessarily have enough space to accommodate many EVs that simultaneously recharge their batteries.

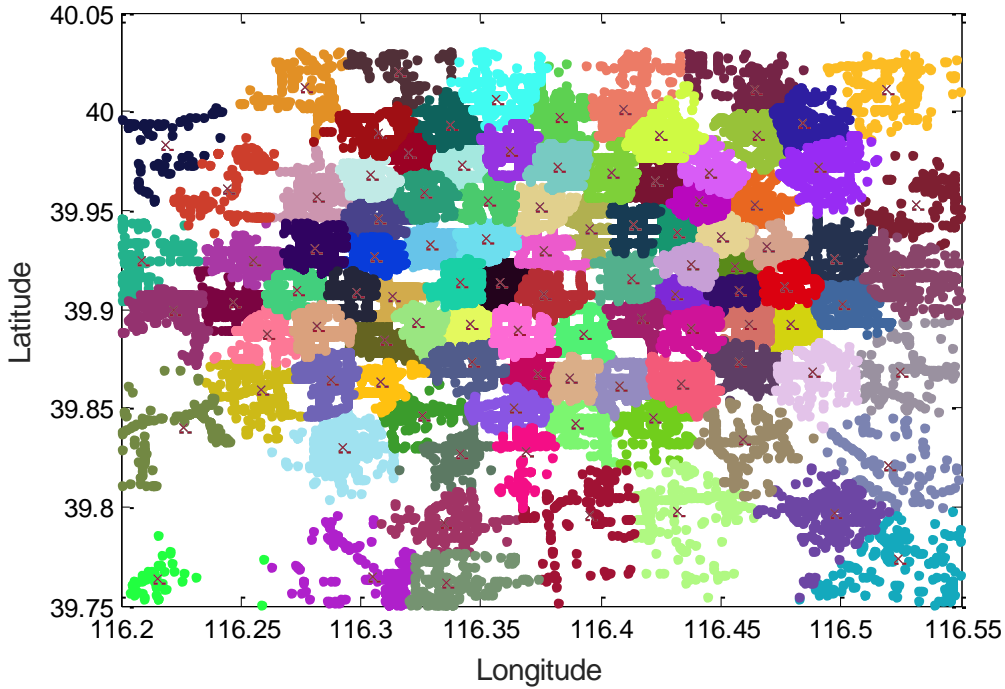
Figure3. The locations of charging stations.

## 3. TIME-SERIES SIMULATIONS AND REAL-TIME FLEET DISPATCHING PLATFORM

On the basis of the gathered real-time vehicle trajectory data of 39,053 taxis in Beijing, we are devoted to proposing a data-driven approach for electrifying the urban taxi fleet. In this section, we conduct the time-series simulations to derive insights on both the configuration of electric taxi fleet and dispatching strategies. Specifically, in the simulations, we relax the unrealistic assumption of the unchanged travel pattern of taxi drivers after shifting to EVs, and develop a centralized fleet-dispatching platform to serve the taxi customer demands, extracted from the dataset, and arrange the recharging of EVs in real time. Note that because of the limited battery range, we assume that EV taxi drivers prefer to rely on the e-hailing or ride-sourcing smartphone-application (such as Uber or Didi) to search for customers, rather than waste battery ranges on cruising on streets. In such context, taxi customers request rides via the mobile application, and upon receiving the ride requests, the centralized platform will match the electric taxis with the taxi customers.

## 3.1 Time-series Simulations

Figure 4 shows the flow chart of the simulation model. In the simulation, time is discretized. As the time step propagates, the O-D taxi customer demands are dynamically fed into the simulation, as per the specification in the dataset. In order to reduce the computational complexity and facilitate the real-time operations of the centralized platform, we divide the whole city into sub-regions. Specifically, considering the operations of electric taxi fleet own a high reliance on the charging infrastructure, we thus use the locations of charging stations for the sub-region dividing, i.e., all the geographic points, sharing the same nearest charging stations, form a sub-region.[3] Figure 5 shows the 100 divided sub-regions based on the deployed charging stations of figure 3. The details of the simulation flowchart in figure 4 are explained in the following eight steps.

**Step1:** The platform receives a new ride request and locates the sub-region it belongs to. Then, the platform selects all the available taxis in the located sub-region as candidates.

**Step2:** Feed all the candidates into reachability-analysis module and identify the candidates whose battery states of charge (SOC) are sufficient to fulfill the customer's travel request.

**Step3:** If there is no candidate passing the reachability test, go to Step 6. Otherwise, go to Step 4.

**Step4:** Feed all the candidates, passing the reachability test, into grading module, and calculate their scores. Dispatch the candidate with the highest score to pick up the customer.

**Step5:** After the assigned taxi reaches the customer's destination, check if its battery SOC is below a pre-determined threshold. If so, send the taxi to the charging station of the current sub-region. Otherwise, flag the taxi as available state, and the taxi will stay at the current place to wait for the next assigned ride request.

**Step6:** Add the demand into the waiting list. The waiting list follows the first-in-first-out (FIFO) rule. Moreover, at each time step, the ride requests in the waiting list own higher priority than the newly-coming requests, implying that the requests in the waiting list should be satisfied first if there are new taxis available.

**Step7:** When the customer's waiting time exceeds a pre-determined waiting-time threshold, the platform will attempt to dispatch the taxis from the adjacent areas. If there does exist an available taxi, which can pass the reachability test, in adjacent areas, then dispatch the taxi to pick up the customer.

---

[3] Doing this can easily exclude the taxis that are far away from the origin of a taxi demand, hence simplifying the taxi-searching procedure of the centralized platform. Furthermore, the sub-regions can also help identify nearby available chargers to satisfy electric taxis' charging demands. All the distances in this paper refer to the Manhattan distance.

**Step8:** When the customer's waiting time exceeds a pre-determined canceling threshold, decline the ride request.

Note that in step5, in charging stations, EVs can charge batteries only if there are still unoccupied chargers left at stations, and otherwise they need to wait in the queue for next available chargers. Doing this makes our simulations capable of accurately modeling the real-time operations of public charging stations and reflecting the interactions of different EVs' charging behaviors. To realize the above consideration, similar to Li et al. (2017), the simulation model utilizes the station operations chart. Specifically, the station operation chart describes each charging station's real-time state, whose rows and columns respectively correspond to time steps and station IDs. Each column of the chart records the number of a station's vacant chargers at different times. We update the chart each time a charger is occupied by an EV.

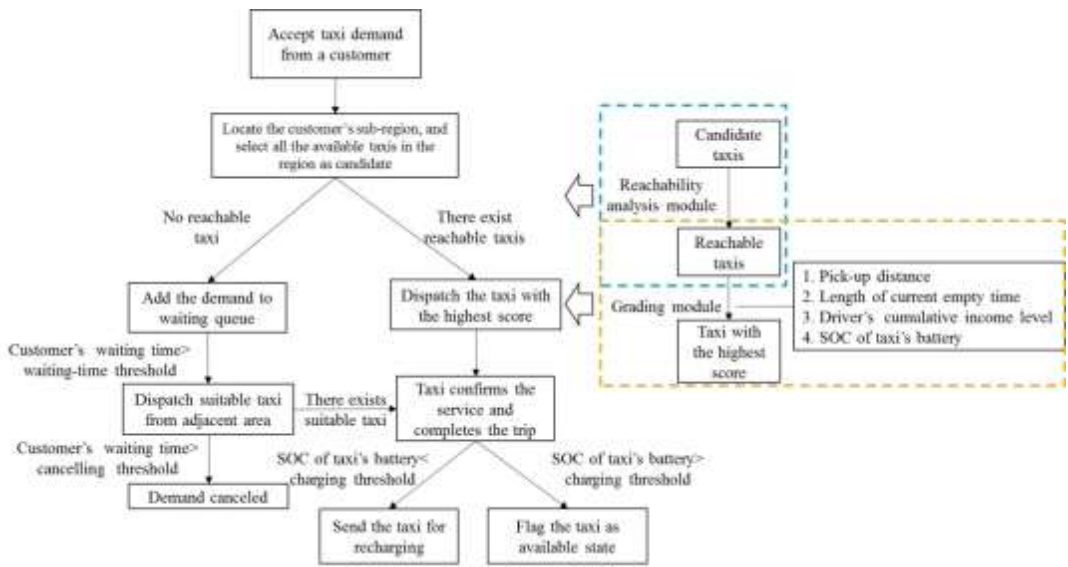

Figure4. The flowchart of the simulation model

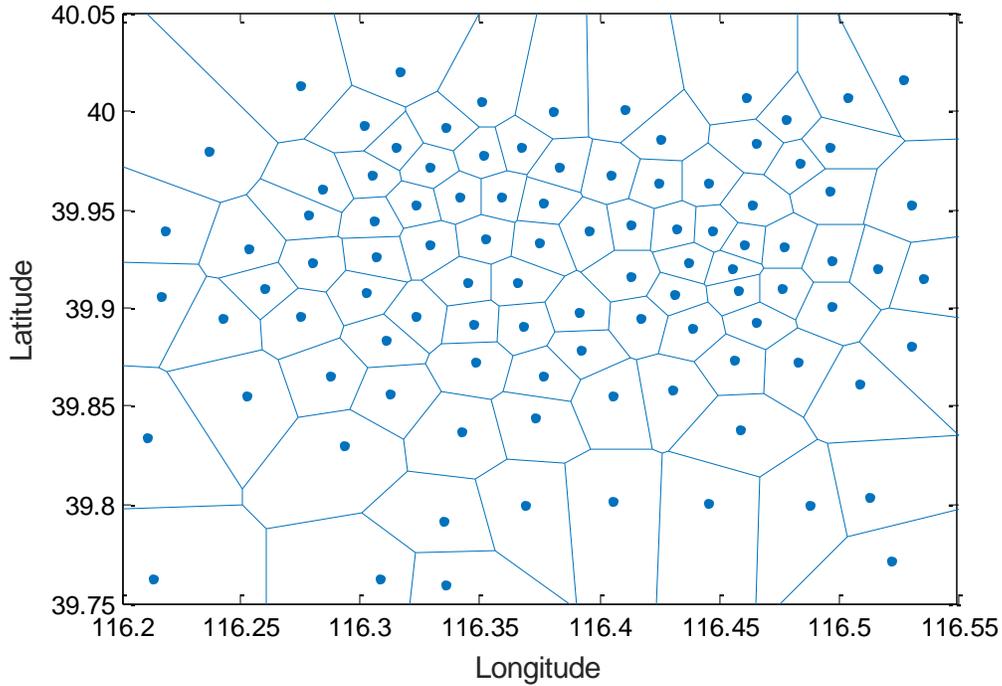

Figure5. The divided sub-regions

**3.2 Rule-based Dispatching Platform.**

To facilitate the real-time operations of the centralized taxi dispatching platform, we adopt a rule-based approach, which considers both efficiency and equity. In terms of efficiency, we aim at reducing customer waiting time and increasing the demand fill rate. With regard to equity, we take into account the equity of taxi drivers' income distribution. The dispatching rule is mainly reflected in the grading module. Recall that in step 4 of the simulation, we feed all the candidates, passing the reachability test, into grading module, calculate their scores and dispatch the candidate with the highest score. Considering both efficiency and equity, the calculation of the score for a taxi in the grading module utilizes the pick-up distance, the length of the taxi's current empty time, the taxi driver's cumulative income level and the taxi's SOC. Specifically, the pick-up distance means the distance between the taxi's current location and the trip's origin. We prefer to dispatch the taxi closer to the customer to reduce the customer's waiting time. The length of the taxi's current empty time measures how long the taxi driver has waited since the finish of his/her last trip. Incorporating it into the grading module implies the platform's consideration of the taxi drivers' willingness to start next trips as early as possible. The taxi driver's cumulative income level reflects the taxi driver's total revenue since an initial time point[4]. Considering equity, the dispatching platform tends to let the driver, with a lower level of cumulative income, own a higher possibility to receive an order of travel. With regard to the SOC of the taxi, we adopt a similar strategy to Lu et al. (2012). If the charging station in the sub-region of the trip's destination is relatively busy, the dispatching platform tends to dispatch a taxi with a higher SOC, in order to avoid possibly long waiting time of recharging. Otherwise, the platform chooses to dispatch a taxi with a lower SOC, for the purpose of better utilizing the charging resource.

Let $E_a$ denote the set of candidate taxis in sub-region $a$, passing the reachability-analysis test. For $e \in E_a$, denote $l_e$, $z_e$ and $o_e$ as the length of taxi's current empty time, cumulative income level and SOC, respectively. For taxi $e$ and trip $n$, if the charging station utilization level in the sub-region of the trip $n$'s destination is above a threshold, we calculate the score through equation (1). Otherwise, we adopt equation (2).

$$S_{e,n} = -w_1 q_1 \hat{d}_{en} + w_2 q_2 l_e - w_3 q_3 \frac{z_e}{T_e} + w_4 q_4 o_e \qquad \forall e, n \qquad (1)$$
$$S_{e,n} = -w_1 q_1 \hat{d}_{en} + w_2 q_2 l_e - w_3 q_3 \frac{z_e}{T_e} - w_4 q_4 o_e \qquad \forall e, n \qquad (2)$$

where $\hat{d}_{en}$ represents the distance between the location of taxi $e$ and the origin of trip $n$; $q_1, q_2, q_3$ and $q_4$ are the parameters used to balance the magnitude of $\hat{d}_{en}, l_e, z_e$, and $o_e$ and make them comparable; $w_1, w_2, w_3$, and $w_4$ denote the weighting parameters associated with the four different factors; $T_e$ represents the cumulative operations time for taxi $e$. Note that adopting different values of $w_1, w_2, w_3$, and $w_4$ corresponds to different dispatching strategies. In section 4, based on the dataset of 39,053 taxis in Beijing, different simulation results are analyzed to derive the insights on the impacts of the values of the weighting parameters.

Step 2 of the simulation feeds all the candidates into reachability-analysis module and identifies the candidates whose battery states of charge (SOC) are sufficient to fulfill the customer's travel request. In the reachability-analysis, we check if a taxi's battery has enough energy to fulfill a trip. For taxi $e$ and trip $n$, if the equations (3)-(4) are satisfied, we conclude that the taxi passes the reachability test.

---

[4] Normally, the initial time point refers to the start of the simulation process.

$$S_e - \hat{d}_{en} - d_n \geq 0 \qquad \forall e, n \qquad (3)$$
$$S_e - \hat{d}_{en} - d_n - \tilde{d}_n \geq 0, \text{if } S_e - \hat{d}_{en} - d_n \leq C \qquad \forall e, n \qquad (4)$$

where $S_e$ represents taxi $e$'s current state of charge; $d_n$ is the travel distance of trip $n$; $\tilde{d}_n$ represents the distance between the destination of trip $n$ and its nearest charging station; $C$ is a pre-determined charging threshold. After dropping off the customer, the taxi will be sent to recharge if its SOC falls below $C$. In the above equations, equation (3) guarantees that the vehicle will not run out of battery before reaching the customer's destination. Equation (4) demonstrates that if the vehicle needs to recharge after dropping off the customer, its remaining electricity should also be enough for reaching the nearest charging station.

### 3.3 Computational Efficiency Analysis

As aforementioned, we divide the region into 100 sub-regions. In this subsection, we provide further analysis on how this division contributes to computational efficiency. In our simulation, we assume that all EVs in a sub-region are coordinated by a local control center. The local center dispatches EVs in its sub-region to meet customers' travel demands, and it will seek for help from adjacent centers when EVs in the area fail to satisfy local demands.

When a new customer demand emerges, the platform will assign the dispatching task to a local control center according to which sub-region the travel demand's origin is located. Then, the local control center will apply the rule-based approach, as described in section 3.2, to evaluate the suitability of each EV in the sub-region. Compared to dispatching all EVs through a centralized control center, i.e., all the EVs are scheduled by one operator with a high dimension of decision variables, the distributed control strategy greatly reduces the computational burden by dividing the onerous dispatching tasks and assigning them to local control centers. As shown in figure 5, we use the Voronoi diagram to partition the entire region, which makes sure that the relatively near EVs will be considered first. To summarize, such distributed and coordinated control strategy enhances the computational efficiency by filtering the customers' travel demands based on their location before making the dispatching strategy of each individual taxi.

## 4. SIMULATION RESULTS

As previously mentioned, we will evaluate the performance of the electrified urban taxi fleet from the aspects of both efficiency and equity. Specifically, we adopt the average customer waiting time and demand fill rate to evaluate efficiency, and use the Gini coefficient of EV taxi drivers' income to evaluate equity (Sen,1973). In this section, different simulation results are compared and analyzed to quantify the impacts of EV number, battery range, station capacity and dispatching-strategy parameters on the performances of electrified taxi fleet.

### 4.1 Simulation Environment

The geographical range of simulation region is from 115.5°E to 117.37°E for longitude and from 39.47°N to 40.68°N for latitude, which forms a 165km ×138km rectangle area. At the beginning of the simulation, the number of taxis in each sub-region is determined by the spatial distribution of taxi demands extracted from the taxi real-time trajectory dataset, and the specific location of each taxi within the sub-region is randomly generated. The taxi fare is calculated as per the current taxi fare standard in Beijing (He et al., 2018). To conduct the sensitivity analyses of the dispatching-strategy parameters, we run the time-series simulations with 16 different values of vector $(w_1, w_2, w_3, w_4)$, as specified in table 1. Note that strategy 16 represents randomly selecting a reachable taxi from the corresponding sub-region to serve the ride request. Lastly, it is assumed that each taxi is operated by two drivers per day, each of which works a 12-hour shift. Table 2 lists the values of other parameters in the simulation.

Table1. 16 different dispatching strategies

| Strategy number | $w_1$ | $w_2$ | $w_3$ | $w_4$ |
|---|---|---|---|---|
| 1 | 1 | 0 | 0 | 0 |
| 2 | 0 | 1 | 0 | 0 |
| 3 | 0 | 0 | 1 | 0 |
| 4 | 0 | 0 | 0 | 1 |
| 5 | 1 | 1 | 0 | 0 |
| 6 | 1 | 0 | 1 | 0 |
| 7 | 1 | 0 | 0 | 1 |
| 8 | 0 | 1 | 1 | 0 |
| 9 | 0 | 1 | 0 | 1 |
| 10 | 0 | 0 | 1 | 1 |
| 11 | 1 | 1 | 1 | 0 |
| 12 | 1 | 1 | 0 | 1 |
| 13 | 1 | 0 | 1 | 1 |
| 14 | 0 | 1 | 1 | 1 |
| 15 | 1 | 1 | 1 | 1 |
| 16 | 0 | 0 | 0 | 0 |

Table2. Parameters of the simulation

| Parameter | Value |
|---|---|
| Taxi speed without passengers | 30km/h |
| Taxi speed with passengers | Calculated based on the dataset |
| Waiting-time threshold | 3min |
| Canceling threshold | 15min[5] |
| Recharging threshold[6] | 20km |
| Recharging time | 30min |
| The threshold for busy stations[7] | 50% |

**4.2 Electric Taxi Fleet Configuration**

In this subsection, we conduct the sensitivity analyses with respect to the number of electric taxis, charging station capacity and battery range. Specifically, we run multiple simulations under various taxi fleet configurations, as shown in table 3. For simplicity, we do not present the detailed results, but summarize the major insights. It is observed that the taxi customer average waiting time (demand fill rate) decreases (increases) with the number of electric taxis, charging station capacity and battery range. Furthermore, we also observe that when EV taxi number is more than 12000, station capacity is greater than 16, and the battery range is more than 200km, there is a marginally diminishing effect. Under the configuration of 12,000 EVs, 16 chargers per station and 200km battery, the customer average waiting time stays below five minutes and the demand fill rate is nearly 100% in most cases.

---

[5] The canceling threshold is set to be 15min unless otherwise stated.
[6] If the EV's SOC falls below 20km, the dispatching platform will send it to the nearest station to charge.
[7] If half of a station's chargers are occupied, we flag it as a busy station. Equation (3) will be adopted to calculate the score.

Table3. Sensitivity analysis of EV fleet configuration

| Parameter | Value |
| --- | --- |
| Time period | May 1$^{st}$ to May 7$^{th}$ of 2016 |
| BEV number | (9000,10000,11000,12000) |
| Range | (100km,150km,200km,250km 300km) |
| Station capacity | (8,16,24,32,40) |
| The number of adjacent areas | 3 |
| Strategy number | (1,2,3,4,6) |

Under dispatching strategy 1, figure 7 depicts the relationship between the Gini coefficient[8] and the battery range under various fleet configurations. In all cases, the Gini coefficient decreases with the battery range. The possible explanation is that the increase of battery range could make the taxis' available time longer, thus providing the dispatching platform more opportunities to balance their incomes. This result suggests that for the taxi fleet in a large city, choosing an enough battery range could not only enhance the efficiency of system but also promote the fairness of drivers' income.

---

[8] We assume that for the two drivers of a taxi, their day or night shift numbers as well as incomes tend to be identical in a long operations term. Thus, we calculate the Gini coefficients for taxis here.

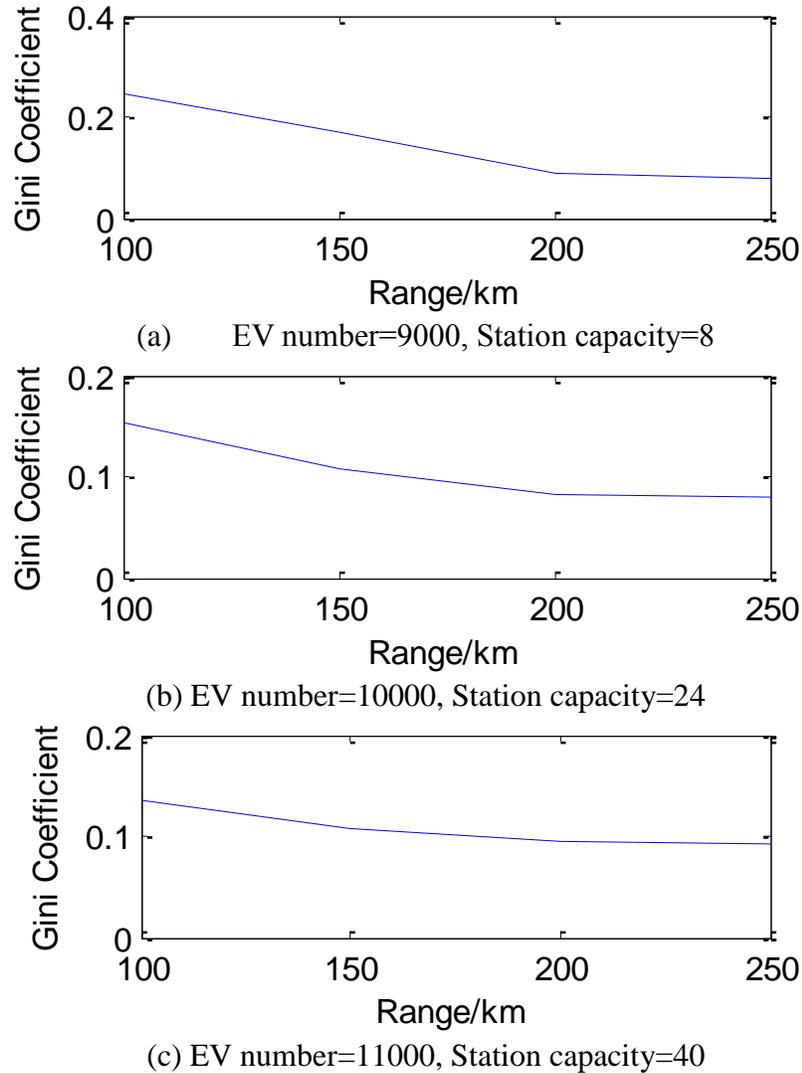

(a) EV number=9000, Station capacity=8

(b) EV number=10000, Station capacity=24

(c) EV number=11000, Station capacity=40

Figure7. The Gini coefficient under various fleet configurations

## 4.3 Dispatching-strategy Parameters

### 4.3.1. Comparisons of different strategies.

We run the simulations with strategies 1-15, under the parameters as specified in table 4. The corresponding average customer waiting times and Gini coefficients are compared in figure 8[9]. The dispatching strategies including the non-zero weighting factor for pick-up distance, i.e., $w_1$, all correspond to a customers' average waiting time of around two minutes or less, implying that considering pick-up distance definitely contributes to the system efficiency. Moreover, among the first eight strategies, (1; 1; 1; 0), (1; 0; 1; 0), (1; 1; 1; 1) and (1; 0; 1; 1) successfully show their capabilities of keeping relatively low levels of Gini coefficients. These four strategies all consider the drivers' cumulative income levels, i.e., the value of $w_3$ is one. In summary, considering both efficiency and equity, dispatching strategies (1; 1; 1; 0), (1; 0; 1; 0), (1; 1; 1; 1) and (1; 0; 1; 1) are recommended for this simulation setting.

---

[9] Under the setting of table 4, the system can keep demand fill rate higher than 97%. In such context, much attention should be paid to the customer average waiting time and Gini coefficient.

Table4. Simulation parameters

| Parameter | Value |
|---|---|
| Time period | May 1st to May 7th of 2016 |
| BEV number | 12000 |
| Range | 200km |
| Station capacity | 16 |
| The number of adjacent areas | 3 |
| Strategy number | (1,2,3,...,15) |

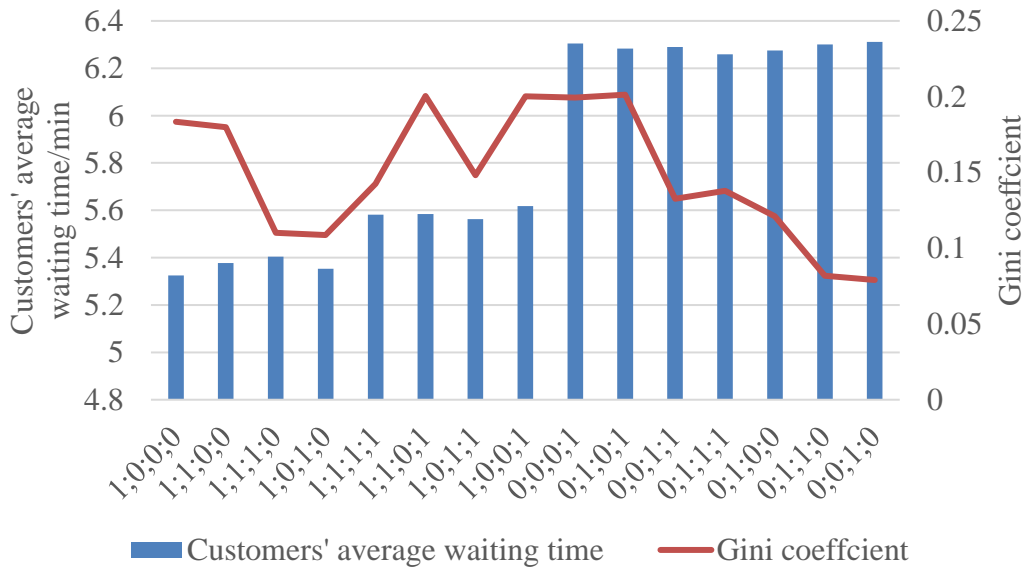

Figure8. Customer average waiting times and Gini coefficients

We are also interested in the performances of those strategies when the charging station capacity and EV number are not so sufficient as in table 4. Setting EV number as 10,000, station capacity as eight and battery range as 150km, table 5 compares the performance of four simple strategies. Compared to the other two strategies, strategies (1; 0; 0; 0) and (0; 0; 0; 1) not only perform better in customer average waiting time and demand fill rate, but also promote the fairness of drivers' income, suggesting that the dispatching strategy should dynamically fit the electric taxi fleet configuration. When the EV number and station capacity are at low levels, priority should be given to the operations efficiency in designing dispatching-strategy parameters. Moreover, the effect of only considering equity is limited when the system efficiency fails to be maintained.

Table 5. The performance of four different dispatching strategies

| Strategy | Gini coefficient | Customer average waiting time (min) | Demand fill rate |
|---|---|---|---|
| 1;0;0;0 | 0.157 | 8.21 | 93.64% |
| 0;1;0;0 | 0.167 | 9.30 | 91.91% |
| 0;0;1;0 | 0.158 | 9.48 | 90.82% |
| 0;0;0;1 | 0.134 | 7.33 | 97.63% |

**4.3.2. Dispatching strategy validation**

To validate the four recommended strategies, i.e., (1; 0; 1; 0), (1; 1; 1; 0), (1; 0; 1; 1) and (1; 1; 1; 1), we compare them with the strategy of randomly selecting a taxi from all the reachable taxis in the sub-region, i.e., strategy (0; 0; 0; 0). Table 6 specifies the simulation parameters. Figure 9 compares the Gini coefficients, demand fill rates and customer average waiting times. Major observations include: i) random dispatching strategy does not fall behind the recommended strategies in terms of demand fill rate when EV number and charging station capacity are deployed sufficiently. In other words, as long as there are enough EV taxis in the city, the demands will be satisfied well whatever the dispatching strategy is; ii) however, the customer average waiting time and the Gini coefficient vary dramatically. Take strategy (1; 0; 1; 1) as an example. Compared to the random dispatching strategy, it decreases the customer average waiting time by 25% and the Gini coefficient by more than 33%, demonstrating that the recommended four strategies remarkably promote the system operation efficiency and the equity of drivers' income.

Table 6. Simulation parameters

| Parameter | Value |
|---|---|
| Time period | May 1$^{st}$ to May 7$^{th}$ of 2016 |
| BEV number | (9000,10000,11000) |
| Range | 100km |
| Station capacity | 8 |
| The number of Adjacent areas | 3 |
| Strategy number | (6,11,13,15,16) |

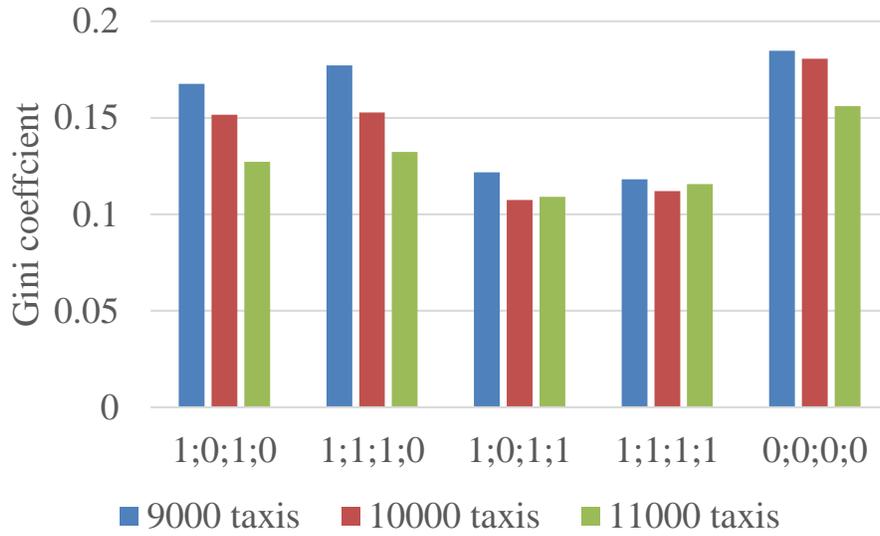

(a) Gini coefficient

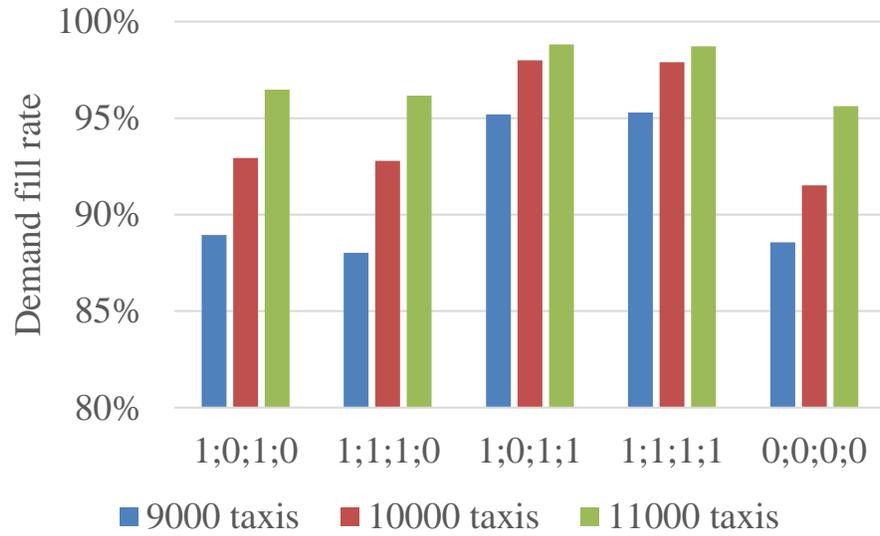

(b) Demand fill rate

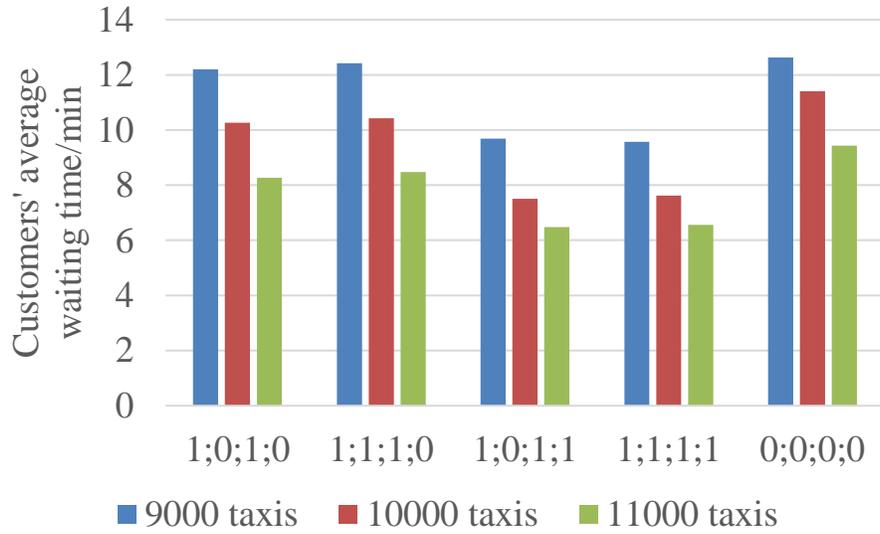

(c) Customer average waiting time

Figure9. Comparisons of different strategies

**4.3.3. The Robustness of dispatching strategies**

It is worthwhile to note that the taxi customer demand density varies with time. In this subsection, we are thereby devoted to investigating the robustness of dispatching strategies under various demand levels of taxi customers. Specifically, based on the dataset, we generate three demand density values (see Table 7) to represent low, middle and high travel demand levels.

Table 7. Simulation parameters

| Parameter | Value |
| --- | --- |
| Demand density | Low demand (1.66 million/week) Middle demand (1.84 million/week) High demand (2.03 million/week) |
| BEV number | 9000 |
| Range | 125km |
| Station capacity | 10 |
| The number of adjacent areas | 3 |
| Strategy number | (1,2,3,…,15) |
| Canceling threshold | 20min |

Figure 10 shows the performances of dispatching strategies under low, middle and high travel demands. It can be observed that the recommended dispatching strategies (1; 0; 1; 0), (1; 1; 1; 0), (1; 0; 1; 1) and (1; 1; 1; 1) all perform well when the trip number is relatively low. However, as the trip number increases, strategies (1;0;1;1) and (1; 1; 1; 1) are still capable of maintaining a good performance, whereas strategies (1; 0; 1; 0) and (1; 1; 1; 0) encounter substantial increases of Gini coefficient, customers' average waiting time and unsatisfied demand rate. Furthermore, figure 10 shows that strategies (0; 1; 0; 1), (0; 0; 0; 1), (1; 1; 0; 1) and (1; 0; 0; 1) also perform well when the travel demand is high. These strategies' common feature is considering SOC of EVs and ignoring the equity of drivers' income. These results further confirm the previous insights that the dispatching strategy should dynamically fit the electric taxi fleet configuration as well as the travel demand level. When demand level is high, priority should be given to concerning the demand satisfaction rate and customer average waiting time. Because of the sufficient taxi customer demands, the equity issue of drivers' income is not significant if the system operations efficiency can be guaranteed.

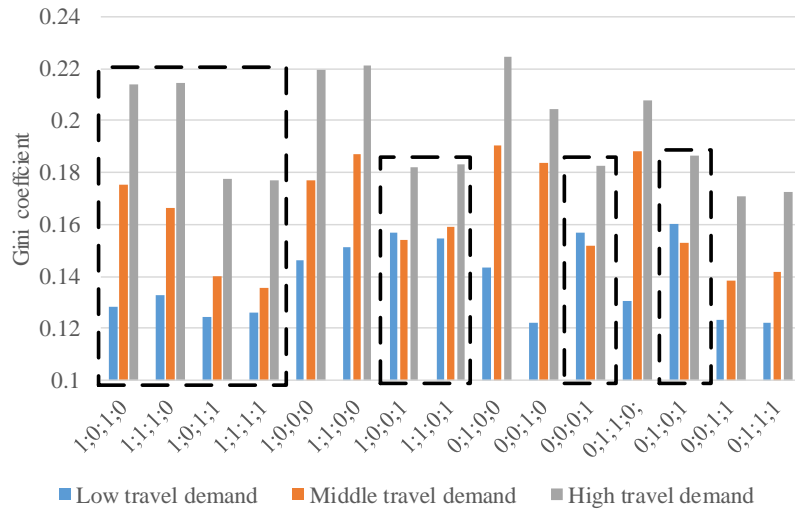

(a) Gini coefficient

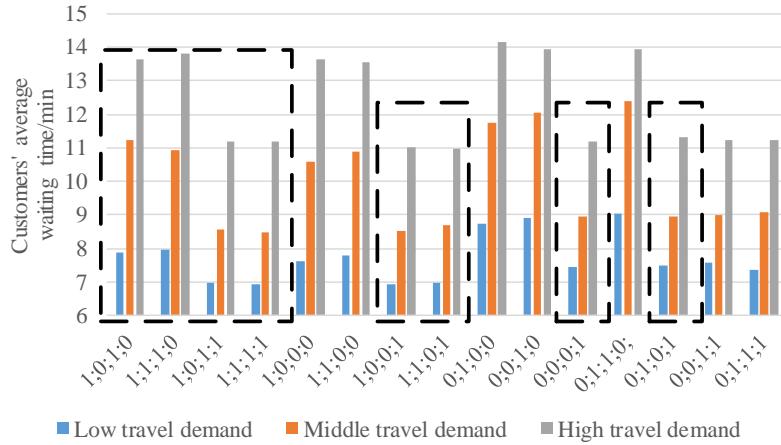

(b) Customer average waiting time

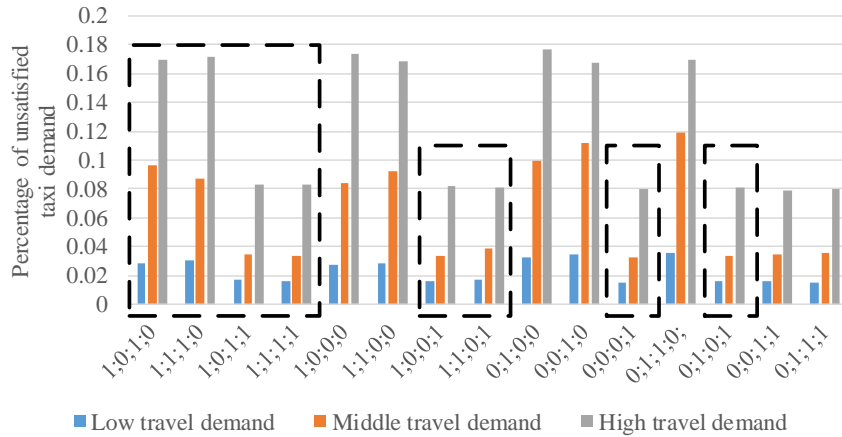

(c) Unsatisfied demand rate

Figure10. The performances of dispatching strategies under low, middle and high customer demands

## 4.4. Temporal Distributions of Taxi Customer Demands and EV Charging Demands

As both the extracted O-D taxi customer demand and the simulation model are dynamic, it would be insightful to see the results in the temporal dimension. For instance, we explore into the relationship between the temporal distributions of customer demands and charging demands. In order to do so, we run a simulation with the parameters listed in table 8. Figure 11 shows the curves of customer and charging demands in the first 1500-minute period of the simulation.

Table8. Simulation parameters

| Parameter | Value |
| --- | --- |
| Time period | May $1^{st}$ to May $7^{th}$ of 2016 |
| BEV number | 12,000 |
| Range | 200km |
| Station capacity | 16 |
| The number of adjacent areas | 3 |
| Strategy number | 6 |

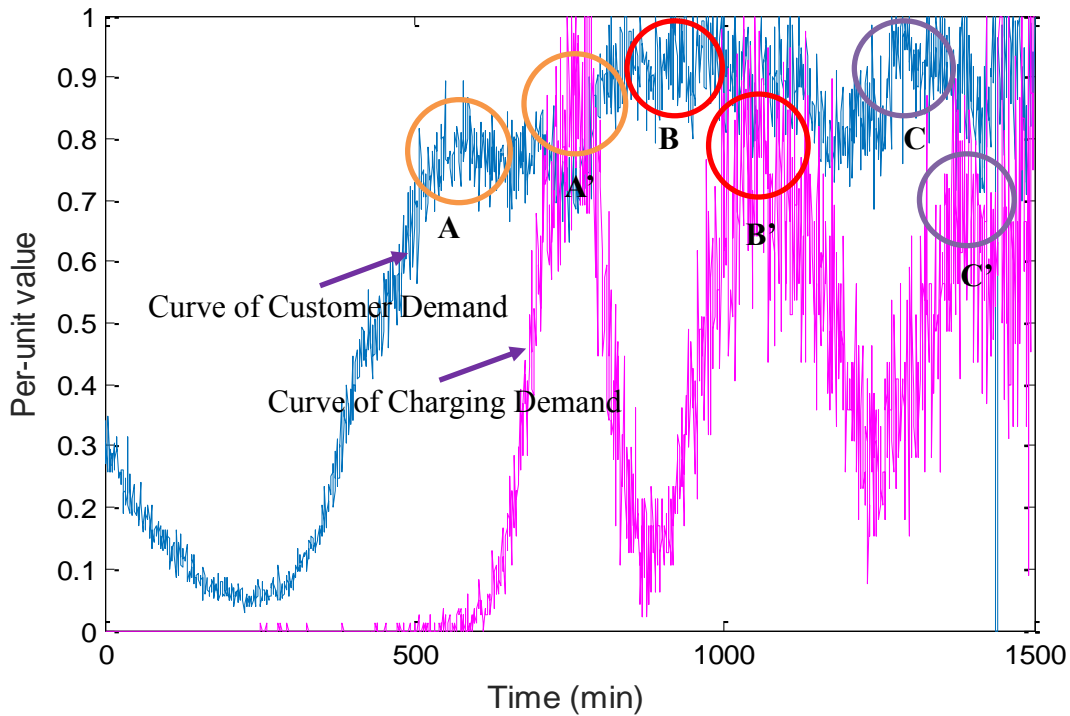

Figure11. The comparison of charging-demand curve and customer-demand curve

In figure 11, the peaks on the curve of customer demand are circled and are labeled as A, B and C, whereas the peaks on the curve of charging demand are circled and labeled as A', B' and C'. It can be observed that peaks A', B' and C' arrive later than peaks A, B and C, respectively, which is reasonable because an EV is only subject to low SOC after delivering several customers. Moreover, according to figure 11, a peak of customer demand is always corresponding to a valley in the charging demand, which could be explained because EV taxis will not consider recharging during service trips. The curve of charging demand shows a much more fluctuated pattern than the curve of customer demand. This phenomenon is likely to impose a significant pressure on the power network, which demands a future research on how to smooth the curve of charging demand so as to mitigate the impact of charging load on the power grid.

## 5. CONCLUSIONS

In this paper, we are devoted to proposing a data-driven approach for electrifying the urban taxi fleet. In order to characterize the spatial and temporal distribution of customer demands, this paper gathers the real-time vehicle trajectory data of 39,053 taxis in Beijing and extracts the customer origin-destination (O-D) pairs whose number reaches approximately eight million in one month. On the basis of this dataset, we conduct time-series simulations to derive insights on both the configuration of electric taxi fleet and dispatching strategies. Specifically, the time-series simulation models the electric vehicle recharging behaviors from the aspects of time window, charging demand and availability of unoccupied charges. Moreover, a centralized fleet-dispatching platform, similar to Uber or Didi, is developed and considered to handle taxi service requests and arrange electric taxis' recharging in real time. To address the impacts of the limited driving range and long battery-recharging time on the electrified fleet's operations efficiency, the dispatching platform integrates the information of customers, taxi drivers and charging stations, and adopts a rule-based approach to achieve multiple objectives, including reducing taxi customer average waiting time, increasing the taxi demand fill rate and guaranteeing the equity of the income among taxi drivers. Simulation results are compared and the major insights are summarized as below.

- The taxi customer average waiting time (demand fill rate) decreases (increases) with the

number of electric taxis, charging station capacity and battery range. On the basis of the gathered vehicle trajectory dataset, when the EV taxi number is more than 12,000, station capacity is greater than 16, and the battery range is more than 200km, continuously increasing them brings marginally diminishing effect for electrifying the taxi fleet in Beijing.
- For the electrified taxi fleet in Beijing, choosing a sufficient battery range could not only enhance the system efficiency but also promote the fairness of drivers' income.
- Under the low and middle demand levels, the dispatching strategies, considering both pick-up distance and drivers' cumulative income level, perform well in terms of efficiency and equity.
- The dispatching strategy should dynamically fit the electric taxi fleet configuration as well as travel demand level. When the EV number and station capacity are at low levels, or the demand level is high, priority should be given to maintaining demand fill rate and reducing customers' average waiting time. Otherwise, not only the efficiency but also equity will be compromised.
- The curve of charging demand shows a much more fluctuated pattern than the curve of customer demand. This phenomenon is likely to impose a significant pressure on the power network.

Future research may consider the impacts of city road network topology and traffic congestion on the time taxis need to pick up customers. Furthermore, this study ignores the impact of EV taxi speed variation on energy consumption rate. We will further extend the simulation framework by adopting more sophisticated models to track SOC of EVs (e.g., Yang et al., 2015). Finally, to better facilitate the urban taxi fleet electrification, future research may also investigate on how to smooth the curve of EV taxi charging demand so as to mitigate the impact of charging loads on the power grid.

## ACKNOWLEDGEMENTS

The research is partially supported by grants from National Natural Science Foundation of China (71501107) and Beijing Natural Science Foundation (8173056). The research is supported in part by the Center for Data-Centric Management in the Department of Industrial Engineering at Tsinghua University.